\documentclass[12p]{amsart}
\usepackage{amssymb}
\usepackage{amsmath}
\usepackage{amsfonts}
\usepackage{geometry}
\usepackage{graphicx}
\usepackage{mathrsfs,amssymb}

\usepackage{hyperref}
\usepackage{cleveref}

\theoremstyle{plain}

\newtheorem{lemma}{Lemma}

\newtheorem{theorem}{Theorem}
\numberwithin{equation}{section}

\begin{document}

\title{Global well-posedness for the defocusing, cubic nonlinear Schr{\"o}dinger equation with initial data in a critical space}

\author{Benjamin Dodson}

\begin{abstract}
In this note we prove global well-posedness for the defocusing, cubic nonlinear Schr{\"o}dinger equation with initial data lying in a critical Sobolev space. 
\end{abstract}
\maketitle
\tableofcontents

\maketitle

\section{Introduction}
In this note, we discuss the defocusing, cubic, nonlinear Schr{\"o}dinger equation in three dimensions,
\begin{equation}\label{1.1}
i u_{t} + \Delta u = F(u) = |u|^{2} u, \qquad u(0,x) = u_{0} \in \dot{H}^{1/2}(\mathbb{R}^{3}).
\end{equation}
Equation $(\ref{1.1})$ has a scaling symmetry. For any $\lambda > 0$, if $u$ solves $(\ref{1.1})$, then
\begin{equation}\label{1.2}
u_{\lambda}(t,x) = \lambda u(\lambda^{2} t, \lambda x),
\end{equation}
also solves $(\ref{1.1})$. The initial data $\lambda u_{0}(\lambda x)$ has $\dot{H}^{1/2}(\mathbb{R}^{3})$ norm that is invariant under the scaling $(\ref{1.2})$.\medskip

The local theory for initial data lying in $\dot{H}^{1/2}(\mathbb{R}^{3})$ has been completely worked out, and the scaling symmetry has been shown to control the local well-posedness theory.
\begin{theorem}\label{t1.1}
Assume $u_{0} \in \dot{H}^{1/2}(\mathbf{R}^{3})$, $\| u_{0} \|_{\dot{H}^{1/2}(\mathbb{R}^{3})} \leq A$. Then there exists $\delta = \delta(A)$ such that if $\| e^{i t \Delta} u_{0} \|_{L_{t,x}^{5}(I \times \mathbb{R}^{3})} < \delta$, then there exists a unique solution to $(\ref{1.1})$ on $I \times \mathbb{R}^{3}$ with $u \in C(I; \dot{H}^{1/2}(\mathbb{R}^{3}))$, and
\begin{equation}\label{1.3}
\| u \|_{L_{t,x}^{5}(I \times \mathbb{R}^{3})} \leq 2 \delta.
\end{equation}
Moreover, if $u_{0, k} \rightarrow u_{0}$ in $\dot{H}^{1/2}(\mathbb{R}^{3})$, the corresponding solutions $u_{k} \rightarrow u$ in $C(I ; \dot{H}^{1/2}(\mathbb{R}^{3}))$.
\end{theorem}
This theorem was proved in \cite{cazenave1990cauchy}.\medskip

From this, it is straightforward to show that local well-posedness holds for $(\ref{1.1})$ for any initial data $u_{0} \in \dot{H}^{1/2}(\mathbb{R}^{3})$. Indeed, by the dominated convergence principle combined with Strichartz estimates, for any $u_{0} \in \dot{H}^{1/2}(\mathbb{R}^{3})$,
\begin{equation}\label{1.4}
\lim_{T \searrow 0} \| e^{it \Delta} u_{0} \|_{L_{t,x}^{5}([-T, T] \times \mathbb{R}^{3})} = 0.
\end{equation}
Since $\delta(A)$ is decreasing as $A \nearrow +\infty$, Strichartz estimates imply that there exists $\delta_{0} > 0$ such that if $\| u_{0} \|_{\dot{H}^{1/2}(\mathbb{R}^{3})} < \delta_{0}$, $(\ref{1.1})$ has a global solution that scatters. By scattering, we mean that there exist $u_{0}^{+}$, $u_{0}^{-}$ so that
\begin{equation}\label{1.5}
\lim_{t \rightarrow +\infty} \| u(t) - e^{it \Delta} u_{0}^{+} \|_{\dot{H}^{1/2}} = 0,
\end{equation}
and
\begin{equation}\label{1.5.1}
\lim_{t \rightarrow -\infty} \| u(t) - e^{it \Delta} u_{0}^{-} \|_{\dot{H}^{1/2}} = 0.
\end{equation}

However, it is important to note that while $(\ref{1.4})$ holds for any fixed $u_{0} \in \dot{H}^{1/2}(\mathbb{R}^{3})$, the convergence is not uniform, even for $\| u_{0} \|_{\dot{H}^{1/2}(\mathbb{R}^{3})} \leq A < \infty$. Thus, one cannot conclude directly from \cite{cazenave1990cauchy} that a uniform bound for $\| u(t) \|_{\dot{H}^{1/2}(\mathbb{R}^{3})}$ on the entire time of the existence of the solution to $(\ref{1.1})$ implies that the solution is global. This result was instead proved in \cite{kenig2010scattering}, using concentration compactness methods.
\begin{theorem}\label{t1.2}
Suppose that $u$ is a solution of $(\ref{1.1})$ with initial data $u_{0} \in \dot{H}^{1/2}(\mathbb{R}^{3})$ and a maximal interval of existence $I = (T_{-}, T_{+})$. Also assume that $\sup_{t \in (T_{-}, T_{+})} \| u(t) \|_{\dot{H}^{1/2}(\mathbb{R}^{3})} = A < \infty$. Then $T_{+}(u_{0}) = +\infty$, $T_{-}(u_{0}) = -\infty$, and the solution $u$ scatters.
\end{theorem}

It is conjectured that $(\ref{1.1})$ is globally well-posed and scattering for any $u_{0} \in \dot{H}^{1/2}(\mathbb{R}^{3})$, without the a priori assumption of a universal bound on the $\dot{H}^{1/2}$ norm of the solution $u(t)$. Partial progress has been made in this direction.

A solution to $(\ref{1.1})$ has the conserved quantities mass,
\begin{equation}\label{1.6}
M(u(t)) = \int |u(t,x)|^{2} dx = M(u(0)),
\end{equation}
and energy,
\begin{equation}\label{1.7}
E(u(t)) = \frac{1}{2} \int |\nabla u(t,x)|^{2} dx + \frac{1}{4} \int |u(t,x)|^{4} dx.
\end{equation}
This fact implies global well-posedness for $(\ref{1.1})$ with $u_{0} \in H_{x}^{1}(\mathbb{R}^{3})$, where $H_{x}^{1}(\mathbb{R}^{3})$ is the inhomogeneous Sobolev space of order one. In this case, one could also prove bounds on the scattering size directly, using the interaction Morawetz estimate of \cite{colliander2004global}.
\begin{theorem}\label{t1.3}
If $u$ is a solution to $(\ref{1.1})$, on an interval $I$, then
\begin{equation}\label{1.8}
\| u \|_{L_{t,x}^{4}(I \times \mathbb{R}^{3})}^{4} \lesssim \| u \|_{L_{t}^{\infty} L_{x}^{2}(I \times \mathbb{R}^{3})}^{2} \| u \|_{L_{t}^{\infty} \dot{H}^{1/2}(I \times \mathbb{R}^{3})}^{2} \lesssim E(u)^{1/2} M(u)^{3/2}.
\end{equation}
\end{theorem}
Interpolating $(\ref{1.7})$ and $(\ref{1.8})$ then implies
\begin{equation}\label{1.9}
\| u \|_{L_{t}^{8} L_{x}^{4}(I \times \mathbb{R}^{3})}^{4} \lesssim M(u)^{3/4} E(u)^{3/4},
\end{equation}
with bounds independent of $I \subset \mathbb{R}$. Combining Strichartz estimates and local well-posedness theory, a uniform bound on $(\ref{1.9})$ for any $I \subset \mathbb{R}$ directly implies a uniform bound on
\begin{equation}\label{1.10}
\| u \|_{L_{t,x}^{5}(I \times \mathbb{R}^{3})}.
\end{equation}
The argument from \cite{cazenave1990cauchy} implies that proving scattering is equivalent to proving
\begin{equation}\label{1.10.1}
\| u \|_{L_{t,x}^{5}(\mathbb{R} \times \mathbb{R}^{3})} < \infty.
\end{equation}
Indeed, assuming $(\ref{1.10.1})$ is true, the interval $\mathbb{R}$ may be partitioned into finitely many pieces $J_{k}$ such that
\begin{equation}\label{1.10.2}
\| u \|_{L_{t,x}^{5}(J_{k} \times \mathbb{R}^{3})} \leq \delta.
\end{equation}
Then iterate the argument over the intervals $J_{k}$, which proves scattering.

This argument also shows that a solution to $(\ref{1.1})$ blowing up at a finite time $T_{0} < \infty$ is equivalent to
\begin{equation}\label{1.10.3}
\| u \|_{L_{t,x}^{5}([0, T_{0}) \times \mathbb{R}^{3})} = \infty.
\end{equation}

\noindent \textbf{Remark:} Prior to \cite{colliander2004global}, \cite{bourgain1998scattering} proved scattering using the standard Morawetz estimate.\medskip

\noindent \textbf{Remark:} See \cite{tao2006nonlinear} for more details on Strichartz estimates.\medskip

Many have attempted to lower the regularity needed in order to prove global well-posedness. For any $s > \frac{1}{2}$, the inhomogeneous Sobolev space $H_{x}^{s}(\mathbb{R}^{3}) \subset \dot{H}^{1/2}(\mathbb{R}^{3})$. Therefore, if $u_{0} \in H_{x}^{s}(\mathbb{R}^{3})$, then it would be conjectured that the solution to $(\ref{1.1})$ with initial data $u_{0}$ is global and scatters.

Proving a uniform bound on the $H_{x}^{s}(\mathbb{R}^{3})$ norm would be enough, since by interpolation this would guarantee a uniform bound on the $\dot{H}_{x}^{1/2}(\mathbb{R}^{3})$ norm. The difficulty is that there does not exist a conserved quantity at regularity $s$ that controls the $\dot{H}^{s}$ norm for $\frac{1}{2} < s < 1$.

Instead, \cite{bourgain1998refinements} used the Fourier truncation method. Decompose the initial data
\begin{equation}\label{1.11}
u_{0} = P_{\leq N} u_{0} + P_{> N} u_{0} = v_{0} + w_{0}.
\end{equation}
Then $v_{0} \in H^{1}(\mathbb{R}^{3})$, and $\| w_{0} \|_{\dot{H}^{1/2}(\mathbb{R}^{3})}$ is small. Thus, $(\ref{1.1})$ has a global solution for initial data $v_{0}$ or $w_{0}$, call them $v$ and $w$. Since $(\ref{1.1})$ is a nonlinear equation, it is necessary to also estimate the interaction between $v$ and $w$ in the nonlinearity of $(\ref{1.1})$. Then, \cite{bourgain1998refinements} proved global well-posedness for $(\ref{1.1})$ with initial data $u_{0} \in H_{x}^{s}(\mathbb{R}^{3})$ when $s > \frac{11}{13}$. Moreover, \cite{bourgain1998refinements} proved that the solution is of the form
\begin{equation}\label{1.12}
e^{it \Delta} u_{0} + v(t),
\end{equation}
where $v(t) \in H_{x}^{1}(\mathbb{R}^{3})$.\medskip

The results of \cite{bourgain1998refinements} for $(\ref{1.1})$ were improved using the I-method. First, \cite{colliander2002almost} improved the regularity necessary for global well-posedness to $s > \frac{5}{6}$. Then, \cite{colliander2004global} improved the necessary regularity to $s > \frac{4}{5}$. To the author's best knowledge, the best known regularity result is the result of \cite{su2011global}, proving global well-posedness and scattering for regularity $s > \frac{5}{7}$. For radial initial data, \cite{dodson2014global} proved global well-posedness and scattering for any $s > \frac{1}{2}$. This result is almost sharp.\medskip

In this paper, we study the cubic nonlinear Schr{\"o}dinger equation $(\ref{1.1})$ with initial data lying in the Sobolev space $W_{x}^{\frac{7}{6}, \frac{11}{7}}(\mathbb{R}^{3})$. That is,
\begin{equation}\label{1.12}
\| |\nabla|^{\frac{11}{7}} u_{0} \|_{L^{\frac{7}{6}}(\mathbb{R}^{3})} < \infty.
\end{equation}
\noindent \textbf{Remark:} This norm is well-defined using the Littlewood--Paley decomposition. See for example \cite{taylor2010partial}.\medskip

This norm is preserved under the scaling $(\ref{1.2})$, and is therefore a critical Sobolev norm. Moreover, $W_{x}^{\frac{7}{6}, \frac{11}{7}}(\mathbb{R}^{3}) \subset \dot{H}^{1/2}(\mathbb{R}^{3})$, so $(\ref{1.1})$ has a local solution for this initial data. We prove global well-posedness for $(\ref{1.1})$ with this initial data.
\begin{theorem}\label{t1.1}
The cubic nonlinear Schr{\"o}dinger equation is globally well-posed for initial data $u_{0} \in W_{x}^{\frac{7}{6}, \frac{11}{7}}(\mathbb{R}^{3})$.
\end{theorem}

The proof of this theorem will heavily utilize dispersive estimates. Interpolating between the fact that $e^{it \Delta}$ is a unitary operator,
\begin{equation}\label{1.13}
\| e^{it \Delta} u_{0} \|_{L^{2}(\mathbb{R}^{3})} = \| u_{0} \|_{L^{2}(\mathbb{R}^{3})},
\end{equation}
and the dispersive estimate,
\begin{equation}\label{1.14}
\| e^{it \Delta} u_{0} \|_{L^{\infty}(\mathbb{R}^{3})} \lesssim \frac{1}{t^{3/2}} \| u_{0} \|_{L^{1}(\mathbb{R}^{3})},
\end{equation}
gives the estimate
\begin{equation}\label{1.15}
\| e^{it \Delta} u_{0} \|_{L^{7}(\mathbb{R}^{3})} \lesssim \frac{1}{t^{\frac{15}{14}}} \| u_{0} \|_{L^{\frac{7}{6}}(\mathbb{R}^{3})}.
\end{equation}

This implies that the linear solution $e^{it \Delta} u_{0}$ has very good behavior when $t > 1$, in fact it is integrable if true. We then rescale so that $u_{0}$ has a local solution on an interval $[-1, 1]$. We prove that this solution may be decomposed into
\begin{equation}\label{1.16}
u(t) = e^{it \Delta} u_{0} + v(t) + w(t).
\end{equation}
In particular,
\begin{equation}\label{1.17}
u(1) = e^{i \Delta} u_{0} + v(1) + w(1).
\end{equation}
The term
\begin{equation}\label{1.18}
e^{i(t - 1) \Delta} e^{i \Delta} u_{0} = e^{it \Delta} u_{0}
\end{equation}
has good properties when $t > 1$. We can also show that
\begin{equation}\label{1.19}
\| \nabla e^{i(t - 1) \Delta} v(1) \|_{L^{\infty}} \lesssim \frac{1}{t^{3/2}},
\end{equation}
which also has good properties when $t > 1$. Finally, $w(1) \in H_{x}^{1}$ and has finite energy. Making a Gronwall argument shows that
\begin{equation}\label{1.20}
\| u(t) - e^{it \Delta} u_{0} - e^{i(t - 1) \Delta} v(1) \|_{\dot{H}^{1}},
\end{equation}
is uniformly bounded on $[1, \infty)$. This is enough to give global well-posedness, but not scattering.\medskip

This result could be compared to the result in \cite{dodson2018global} for the nonlinear wave equation. There, the author proved global well-posedness and scattering for the cubic wave equation with initial radial data in the Besov space $B_{1,1}^{2} \times B_{1,1}^{1}$. Here, we do not require radial symmetry, however, we only prove global well-posedness. We are unable to prove scattering at this time due to the lack of a scale invariant conformal symmetry.\medskip

We prove a local well-posedness result in section two, and a global result in section three. This argument could be generalized to many intercritical, defocusing nonlinear Schr{\"o}dinger equations.

\section{Local well-posedness}
The Sobolev embedding theorem implies that $W_{x}^{\frac{7}{6}, \frac{11}{7}}(\mathbb{R}^{3})$ is embedded into $\dot{H}^{1/2}(\mathbb{R}^{3})$. Therefore, $(\ref{1.1})$ is locally well-posed, and there exists some $T(u_{0}) > 0$ such that $(\ref{1.1})$ has a solution on $[-T, T]$ and $\| u \|_{L_{t}^{5}([-T, T] \times \mathbb{R}^{3})} = \epsilon_{0}$, for some $\epsilon_{0}(\| u_{0} \|_{\dot{H}^{1/2}})$ small. After rescaling using $(\ref{1.2})$, suppose
\begin{equation}\label{2.1}
\| u \|_{L_{t,x}^{5}([-1, 1] \times \mathbb{R}^{3})} = \epsilon_{0}.
\end{equation}
Since $(3, \frac{18}{5})$ is an admissible pair, Strichartz estimates imply
\begin{equation}\label{2.3}
\| |\nabla|^{1/2} u \|_{L_{t}^{\infty} L_{x}^{2} \cap L_{t}^{2} L_{x}^{6}([-1,1] \times \mathbb{R}^{3})} \lesssim \| |\nabla|^{1/2} u_{0} \|_{L_{x}^{2}(\mathbb{R}^{3})} + \| |\nabla|^{1/2} u \|_{L_{t}^{3} L_{x}^{\frac{18}{5}}([-1,1] \times \mathbb{R}^{3})} \| u \|_{L_{t,x}^{5}([-1,1] \times \mathbb{R}^{3})}^{2}.
\end{equation}
Therefore, 
\begin{equation}\label{2.3.1}
\| |\nabla|^{1/2} u \|_{L_{t}^{\infty} L_{x}^{2} \cap L_{t}^{2} L_{x}^{6}([-1,1] \times \mathbb{R}^{3})} \lesssim \| u_{0} \|_{\dot{H}^{1/2}}.
\end{equation}

Also, by Duhamel's principle, for any $t \in [-1, 1]$,
\begin{equation}\label{2.2}
u(t) = e^{it \Delta} u_{0} - i \int_{0}^{t} e^{i(t - \tau) \Delta} F(u(\tau)) d\tau = u_{l}(t) + u_{nl}(t).
\end{equation}

We begin with a technical lemma. This lemma allows us to make a Littlewood--Paley decomposition of $u_{nl}$, treat each $P_{j} u_{nl}$ separately, and then sum up. It also implies that $u_{nl}$ retains all the properties of a solution to the linear Schr{\"o}dinger equation with initial data in a Besov space.\medskip

\noindent \textbf{Remark:} In this section, all implicit constants depend on the norm $\| u_{0} \|_{W^{\frac{7}{6}, \frac{11}{7}}}$.\medskip

\noindent \textbf{Remark:} Throughout this section we rely very heavily on the bilinear Strichartz estimate
\begin{equation}
\| (e^{it \Delta} P_{j} u_{0})(e^{it \Delta} P_{k} v_{0}) \|_{L_{t,x}^{2}(\mathbb{R} \times \mathbb{R}^{3})} \lesssim 2^{-j/2} 2^{k} \| P_{j} u_{0} \|_{L^{2}} \| P_{k} v_{0} \|_{L^{2}}.
\end{equation}
See \cite{bourgain1998refinements} for a proof.
\begin{lemma}\label{l2.1}
Let $P_{j}$ be the customary Littlewood--Paley projection operator. Also suppose that $u$ is a solution to $(\ref{1.1})$ satisfying $(\ref{2.1})$. Then
\begin{equation}\label{2.5}
\sum_{j} 2^{j/2} \| P_{j} F(u) \|_{L_{t}^{1} L_{x}^{2}([-1, 1] \times \mathbb{R}^{3})} \lesssim 1.
\end{equation}
\end{lemma}
\noindent \emph{Proof:} Decompose the nonlinearity,
\begin{equation}\label{2.6}
P_{j} F(u) = P_{j} F(P_{\geq j - 3} u) + 3 P_{j} ((P_{\geq j - 3} u)^{2} (P_{\leq j - 3} u)) + 3 P_{j}((P_{j - 3 \leq \cdot \leq j + 3} u)(P_{\leq j - 3} u)^{2}).
\end{equation}
By Bernstein's inequality, and $(\ref{2.3})$,
\begin{equation}\label{2.7}
\aligned
2^{j/2} \| P_{j} F(P_{\geq j - 3} u) \|_{L_{t}^{1} L_{x}^{2}([-1, 1] \times \mathbb{R}^{3})} \\ \lesssim 2^{j/2} \| P_{\geq j - 3} u \|_{L_{t}^{3} L_{x}^{6}([-1, 1] \times \mathbb{R}^{3})}^{3} \lesssim 2^{j/2} (\sum_{l \geq j - 3} 2^{-l/6} \| |\nabla|^{1/6} P_{l} u \|_{L_{t}^{3} L_{x}^{6}})^{3}.
\endaligned
\end{equation}
Next,
\begin{equation}\label{2.8}
\aligned
2^{j/2} \| P_{j}((P_{\geq j - 3} u)^{2} (P_{\leq j - 3} u)) \|_{L_{t}^{1} L_{x}^{2}([-1, 1] \times \mathbb{R}^{3})}
\lesssim 2^{j/2} (\sum_{l \geq j - 3} 2^{-l/4} \| |\nabla|^{1/4} P_{l} u \|_{L_{t}^{3} L_{x}^{\frac{36}{7}}})^{2} \| u \|_{L_{t}^{3} L_{x}^{9}}.
\endaligned
\end{equation}
Finally, by the bilinear Strichartz estimate and the Sobolev embedding properties of Littlewood--Paley projections,
\begin{equation}\label{2.9}
\aligned
2^{j/2} \| (P_{j - 3 \leq \cdot \leq j + 3} u)(P_{\leq j - 3} u)^{2} \|_{L_{t}^{1} L_{x}^{2}([-1, 1] \times \mathbb{R}^{3})} \\ \lesssim 2^{j/2} \sum_{l_{1} \leq l_{2} \leq j - 3} \| (P_{l_{1}} u)(P_{j - 3 \leq \cdot \leq j + 3} u) \|_{L_{t,x}^{2}} \sum_{l_{1} \leq l_{2} \leq j - 3} \| P_{l_{2}} u \|_{L_{t}^{2} L_{x}^{\infty}} \\
\lesssim 2^{-j/2} \| |\nabla|^{1/2} u \|_{L_{t}^{2} L_{x}^{6}} \sum_{l_{1} \leq j - 3} 2^{l_{1}/2} (j - l_{1}) (\| P_{j - 3 \leq \cdot \leq j + 3} u_{0} \|_{L^{2}} + \| P_{j - 3 \leq \cdot \leq j + 3} F(u) \|_{L_{t}^{1} L_{x}^{2}}) \\ \times (\| P_{l_{1}} u_{0} \|_{L^{2}} + \| P_{l_{1}} F(u) \|_{L_{t}^{1} L_{x}^{2}}).
\endaligned
\end{equation}
By Strichartz estimates, $(\ref{2.3.1})$, Plancherel's theorem, and the fractional product rule,
\begin{equation}\label{2.9.1}
\aligned
\sum_{j} 2^{j} \| P_{j} u_{0} \|_{L^{2}}^{2} + \sum_{j} 2^{j} \| P_{j} F(u) \|_{L_{t}^{1} L_{x}^{2}([-1, 1] \times \mathbb{R}^{3})}^{2} \lesssim \| u_{0} \|_{\dot{H}^{1/2}}^{2} + \| |\nabla|^{1/2} F(u) \|_{L_{t}^{1} L_{x}^{2}}^{2} \\ \lesssim \| u_{0} \|_{\dot{H}^{1/2}}^{2} + \| |\nabla|^{1/2} u \|_{L_{t}^{3} L_{x}^{18/5}}^{2} \| u \|_{L_{t}^{3} L_{x}^{9}}^{4} \lesssim 1.
\endaligned
\end{equation}
Combining $(\ref{2.7})$--$(\ref{2.9})$ with the Cauchy--Schwarz inequality implies
\begin{equation}\label{2.10}
\sum_{j} 2^{j/2} \| P_{j} F(u) \|_{L_{t}^{1} L_{x}^{2}([-1, 1] \times \mathbf{R}^{3})} \lesssim 1,
\end{equation}
which proves the lemma. $\Box$\medskip

Next, decompose $u_{nl}$ in the following manner:
\begin{equation}\label{2.11}
u_{nl}(t) = -i \int_{0}^{(1 - \delta) t} e^{i(t - \tau) \Delta} F(u(\tau)) d\tau - i \int_{(1 - \delta) t}^{t} e^{i(t - \tau) \Delta} F(u(\tau)) d\tau = v(t) + w(t),
\end{equation}
for some $\delta > 0$ sufficiently small, to be specified later.

\begin{lemma}\label{l2.2}
For any $t \in [0, 1]$,
\begin{equation}\label{2.11.1}
\| v(t) \|_{L^{\infty}} \lesssim \frac{1}{\delta^{1/2} t^{1/2}},
\end{equation}
and
\begin{equation}\label{2.12}
\| \nabla v(t) \|_{L^{\infty}} \lesssim \frac{1}{\delta t}.
\end{equation}
\end{lemma}
\noindent \emph{Proof:} By the dispersive estimate, since $\| u \|_{L^{3}} \lesssim \| u \|_{\dot{H}^{1/2}}$ is uniformly bounded on $[0, 1]$,
\begin{equation}\label{2.13}
\| v(t) \|_{L^{\infty}} \lesssim \| \int_{0}^{(1 - \delta) t} e^{i(t - \tau) \Delta} F(u) d\tau \|_{L^{\infty}} \lesssim \int_{0}^{(1 - \delta) t} \frac{1}{|t - \tau|^{3/2}} \| u \|_{L^{3}}^{3} d\tau \lesssim \frac{1}{\delta^{1/2} t^{1/2}}.
\end{equation}

To prove $(\ref{2.12})$, observe that by the product rule,
\begin{equation}\label{2.14}
\nabla F(u) = 2 |u|^{2} \nabla u + u^{2} \nabla \bar{u}.
\end{equation}
Interpolating,
\begin{equation}\label{2.15}
\| |\nabla|^{1/2} u_{l} \|_{L^{2}} \lesssim \| |\nabla|^{1/2} u_{0} \|_{L^{2}} \lesssim 1,
\end{equation}
with
\begin{equation}\label{2.16}
t^{15/14} \| |\nabla|^{11/7} u_{l} \|_{L^{7}} \lesssim \| |\nabla|^{11/7} u_{0} \|_{L^{7/6}} \lesssim 1,
\end{equation}
we have
\begin{equation}\label{2.17}
t^{1/2} \| \nabla u_{l} \|_{L^{3}} \lesssim 1.
\end{equation}
Making a dispersive estimate,
\begin{equation}\label{2.18}
\aligned
\| \int_{0}^{(1 - \delta) t} e^{i(t - \tau) \Delta} |u|^{2} \nabla u_{l}(\tau) d\tau \|_{L^{\infty}} \lesssim \int_{0}^{(1 - \delta) t} \frac{1}{|t - \tau|^{3/2}} \| \nabla u_{l}(\tau) \|_{L^{3}} \| u \|_{L^{3}}^{2} d\tau \\
\lesssim \int_{0}^{(1 - \delta) t} \frac{1}{|t - \tau|^{3/2}} \frac{1}{|\tau|^{1/2}} d\tau \lesssim \frac{1}{\delta t}.
\endaligned
\end{equation}
The same computation may also be made for $u^{2} \nabla \bar{u}_{l}$.\medskip

Next, consider the contribution of $|u|^{2} \nabla u_{nl}$. By $(\ref{2.5})$, we can, without loss of generality, consider only one $P_{j}$ Littlewood-Paley multiplier, provided the estimate is uniform in $2^{j/2} \| P_{j} F(u) \|_{L_{t}^{1} L_{x}^{2}}$.
\begin{equation}\label{2.19}
|u|^{2} (\nabla P_{j} u_{nl}) = |P_{\leq j} u|^{2} (\nabla P_{j} u_{nl}) + 2 Re((P_{> j} \bar{u})(P_{\leq j} \bar{u})) (\nabla P_{j} u_{nl}) + |P_{> j} u|^{2} (\nabla P_{j} u_{nl}).
\end{equation}
Making a bilinear Strichartz estimate and the Cauchy--Schwartz inequality,
\begin{equation}\label{2.20}
\aligned
\| |u_{\leq j}|^{2} (\nabla P_{j} u_{nl}) \|_{L_{t}^{2} L_{x}^{1}([0, 1] \times \mathbf{R}^{3})} \lesssim \sum_{j_{1} \leq j_{2} \leq j} \| (P_{j_{1}} u)(P_{j} \nabla u_{nl}) \|_{L_{t,x}^{2}} \| P_{j_{2}} u \|_{L_{t}^{\infty} L_{x}^{2}} \\
\lesssim \sum_{j_{1} \leq j_{2} \leq j} 2^{j_{1}/2} 2^{-j_{2}/2} 2^{j/2} \| P_{j} F(u) \|_{L_{t}^{1} L_{x}^{2}} (\| |\nabla|^{1/2} P_{j_{1}} u_{0} \|_{L^{2}} + \| |\nabla|^{1/2} P_{j_{1}} F(u) \|_{L_{t}^{1} L_{x}^{2}}) \\ \times (\| |\nabla|^{1/2} P_{j_{2}} u_{0} \|_{L^{2}} + \| |\nabla|^{1/2} P_{j_{2}} F(u) \|_{L_{t}^{1} L_{x}^{2}}) \lesssim 1.
\endaligned
\end{equation}
Also, by Bernstein's inequality,
\begin{equation}\label{2.21}
\aligned
\| |\nabla P_{j} u_{nl}| |P_{> j} u| (|P_{\leq j} u| + |P_{> j} u|) \|_{L_{t}^{2} L_{x}^{1}} \\ \lesssim \| |\nabla|^{1/2} P_{j} u_{nl} \|_{L_{t}^{2} L_{x}^{6}} \| |\nabla|^{1/2} P_{> j} u \|_{L_{t}^{\infty} L_{x}^{2}} \| u \|_{L_{t}^{\infty} L_{x}^{3}} \lesssim 1.
\endaligned
\end{equation}
Therefore,
\begin{equation}\label{2.22}
\aligned
\| \int_{0}^{(1 - \delta) t} e^{i(t - \tau) \Delta} |u|^{2} \nabla u_{nl}(\tau) d\tau \|_{L^{\infty}} \lesssim \int_{0}^{(1 - \delta) t} \frac{1}{|t - \tau|^{3/2}} \| |u|^{2} \nabla u_{nl} \|_{L^{1}} d\tau \\
\lesssim \frac{1}{\delta t} \| |u|^{2} \nabla u_{nl} \|_{L_{t}^{2} L_{x}^{1}} \lesssim \frac{1}{\delta t}.
\endaligned
\end{equation}
The same computation can be also be made for $u^{2} \nabla \bar{u}_{nl}$. This completes the proof of Lemma $\ref{l2.2}$. $\Box$

\begin{lemma}\label{l2.3}
For any $t \in [0, 1]$,
\begin{equation}\label{2.23}
\| |\nabla|^{1/2} w(t) \|_{L^{3}} \lesssim \frac{1}{\delta^{1/4} t^{1/4}}.
\end{equation}
\end{lemma}
\emph{Proof:} First observe that by interpolation, Bernstein's inequality, and $(\ref{2.17})$,
\begin{equation}\label{2.24}
\| |\nabla|^{1/2} e^{it \Delta} u_{0} \|_{L^{3}} \lesssim t^{1/4} \| \nabla e^{it \Delta} P_{\geq t^{-1/2}} u_{0} \|_{L^{3}} + t^{-1/4} \| P_{\leq t^{-1/2}} u_{0} \|_{\dot{H}^{1/2}} \lesssim t^{-1/4}.
\end{equation}
Also since $e^{it \Delta}$ is unitary in $L^{2}$, by $(\ref{2.3})$,
\begin{equation}\label{2.25}
\| v(t) \|_{\dot{H}^{1/2}} = \| u_{nl}((1 - \delta)t) \|_{\dot{H}^{1/2}} \lesssim \| |\nabla|^{1/2} u \|_{L_{t}^{3} L_{x}^{\frac{18}{5}}} \| u \|_{L_{t,x}^{5}}^{2} \lesssim \epsilon_{0}^{2}.
\end{equation}
so interpolating $(\ref{2.11.1})$, $(\ref{2.12})$, and $(\ref{2.25})$,
\begin{equation}\label{2.26}
\| |\nabla|^{1/2} v \|_{L^{3}} \lesssim \| |\nabla|^{1/2} v \|_{L^{\infty}}^{1/3} \| |\nabla|^{1/2} v \|_{L^{2}}^{2/3} \lesssim \frac{\epsilon_{0}^{4/3}}{\delta^{1/4} t^{1/4}}.
\end{equation}
Finally, making a dispersive estimate, for any $t \in [0, 1]$, by $(\ref{2.24})$ and $(\ref{2.26})$, if $\delta^{1/4} \ll \epsilon_{0}$,
\begin{equation}\label{2.27}
\aligned
\delta^{1/4} t^{1/4} \| \int_{(1 - \delta) t}^{t} e^{i(t - \tau) \Delta} |\nabla|^{1/2} F(u) d\tau \|_{L^{3}} \lesssim \delta^{1/4} t^{1/4} \int_{(1 - \delta) t}^{t} \frac{1}{|t - \tau|^{1/2}} \| |\nabla|^{1/2} u(\tau) \|_{L^{3}} \| u(\tau) \|_{L^{6}}^{2} d\tau \\
\lesssim (\sup_{t \in [0, 1]} \delta^{1/4} t^{1/4} \| |\nabla|^{1/2} u \|_{L^{3}})^{3} \lesssim \epsilon_{0}^{4} + (\sup_{t \in [0, 1]} \delta^{1/4} t^{1/4} \| |\nabla|^{1/2} w \|_{L^{3}})^{3}.
\endaligned
\end{equation}
Therefore, absorbing the second term on the right hand side into the left hand side of $(\ref{2.27})$ proves $(\ref{2.23})$.
\begin{equation}\label{2.28}
\| |\nabla|^{1/2} w(t) \|_{L^{3}}  \lesssim \frac{\epsilon_{0}^{4}}{\delta^{1/4} t^{1/4}}.
\end{equation}
 $\Box$\medskip

\noindent \textbf{Remark:} To make the proofs of Lemmas $\ref{l2.2}$ and $\ref{l2.3}$ completely rigorous, truncate $u_{0}$ in frequency. Then the bounds $(\ref{2.11.1})$, $(\ref{2.12})$, and $(\ref{2.23})$ all hold on some open subset of $[0, 1]$ that contains $0$. Making the bootstrap argument using the proof of Lemma $\ref{l2.3}$ gives bounds on all of $[0, 1]$ that do not depend on the frequency truncation of $u_{0}$. Standard perturbation arguments then give the Lemmas.\medskip

Lemma $\ref{l2.3}$ can be strengthened to an estimate on the $\dot{H}^{1}$ norm of $w$.
\begin{lemma}\label{l2.4}
For any $t \in [0, 1]$,
\begin{equation}\label{2.29}
\| \nabla v \|_{L^{2}} \lesssim \frac{1}{\delta^{1/4} t^{1/4}}.
\end{equation}
\end{lemma}
\noindent \emph{Proof:} Once again make use of the bilinear Strichartz estimate. Again by the product rule,
\begin{equation}\label{2.30}
\nabla F(u) = 2 |u|^{2} \nabla u + u^{2} \nabla \bar{u}.
\end{equation}
First, by Strichartz estimates, $(\ref{2.17})$, Lemma $\ref{l2.3}$, and the Sobolev embedding theorem,
\begin{equation}\label{2.31}
\aligned
\| \int_{(1 - \delta)t}^{t} e^{i(t - \tau) \Delta} [2 |u|^{2} \nabla u_{l} + u^{2} \nabla \bar{u}_{l}] d\tau \|_{L^{2}} \lesssim \| 2 |u|^{2} \nabla u_{l} + u^{2} \nabla \bar{u}_{l} \|_{L_{t}^{2} L_{x}^{6/5}} \\ \lesssim \delta^{1/2} t^{1/2} \| \nabla u_{l} \|_{L_{t}^{\infty} L_{x}^{3}([(1 - \delta) t, t] \times \mathbb{R}^{3})} \| u \|_{L_{t}^{\infty} L_{x}^{3}([(1 - \delta) t, t] \times \mathbb{R}^{3})} \| |\nabla|^{1/2} u \|_{L_{t}^{\infty} L_{x}^{3}([(1 - \delta)t, t] \times \mathbb{R}^{3})} \lesssim \frac{\delta^{1/4}}{t^{1/4}}.
\endaligned
\end{equation}
Next, by $(\ref{2.23})$, bilinear Strichartz estimates and the Littlewood-Paley theorem,
\begin{equation}\label{2.32}
\aligned
\| 2 |u_{\leq j}|^{2} (\nabla P_{j} u_{nl}) + (u_{\leq j})^{2} (\nabla P_{j} \bar{u}_{nl}) \|_{L_{t}^{2} L_{x}^{6/5}} \\ \lesssim \sum_{k \geq 0} 2^{-k/2} \| (\sum_{j_{1} \leq j} 2^{j_{1} + k} |P_{j_{1} + k} u|^{2})^{1/2} (\sum_{j_{1} \leq j} 2^{-j_{1}} |P_{j_{1}} u|^{2} |P_{j} u|^{2})^{1/2} \|_{L_{t}^{2} L_{x}^{6/5}} \\
\lesssim \sum_{k \geq 0} 2^{-k/2} \| |\nabla|^{1/2} u(t) \|_{L_{t}^{\infty} L_{x}^{3}([(1 - \delta)t, t] \times \mathbb{R}^{3})} (\sum_{j_{1} \leq j} \| P_{j_{1}} u_{0} \|_{\dot{H}^{1/2}}^{2} + \| P_{j_{1}} F(u) \|_{L_{t}^{1} L_{x}^{2}}^{2})^{1/2} \| P_{j} F(u) \|_{L_{t}^{1} L_{x}^{2}} \\
\lesssim \frac{1}{\delta^{1/4} t^{1/4}} \| |\nabla|^{1/2} P_{j} F(u) \|_{L_{t}^{1} L_{x}^{2}}.
\endaligned
\end{equation}
Next, by Bernstein's inequality,
\begin{equation}\label{2.33}
\aligned
\| (\nabla P_{j} u_{nl}) |u_{\geq j}| |u| \|_{L_{t}^{2} L_{x}^{6/5}} \lesssim \delta^{1/4} t^{1/4} \| |\nabla|^{1/2} u \|_{L_{t}^{\infty} L_{x}^{3}([(1 - \delta) t, t] \times \mathbb{R}^{3})}^{2} \| |\nabla|^{1/2} P_{j} u_{nl} \|_{L_{t}^{4} L_{x}^{3}([(1 - \delta) t, t] \times \mathbb{R}^{3})} \\
\lesssim \frac{1}{\delta^{1/4} t^{1/4}} \| |\nabla|^{1/2} P_{j} F(u) \|_{L_{t}^{1} L_{x}^{2}([0, 1] \times \mathbb{R}^{3})}.
\endaligned
\end{equation}
Summing up in $j$ using Lemma $\ref{l2.1}$ completes the proof. $\Box$\medskip

\noindent \textbf{Remark:} The above arguments would work equally well in the time interval $[-1, 0]$.

\section{Global well-posedness}
We are now ready to prove Theorem $\ref{t1.1}$. The proof will use conservation of the energy $(\ref{1.7})$. Decompose
\begin{equation}\label{3.1}
u(1) = \tilde{v}(1) + w(1),
\end{equation}
where
\begin{equation}\label{3.2}
\tilde{v}(1) = u_{l}(1) + v(1),
\end{equation}
and $w(1)$ is the $w$ in the previous section. Let $T_{0} > 1$ be a time value for which we know that $(\ref{1.1})$ has a solution on $[0, T_{0})$. By standard local well-posedness arguments and 
$(\ref{2.8})$, we know that such a $T_{0}$ exists. Then on $[1, T_{0})$, decompose
\begin{equation}\label{3.3}
u(t) = \tilde{v}(t) + w(t),
\end{equation}
where $\tilde{v}(t)$ is the solution to
\begin{equation}\label{3.4}
(i \partial_{t} + \Delta) \tilde{v}(t) = 0, \qquad \tilde{v}(1) = \tilde{v}(1, x),
\end{equation}
and $w(t)$ is the solution to
\begin{equation}\label{3.5}
(i \partial_{t} + \Delta) w = |u|^{2} u, \qquad w(1) = w(1, x).
\end{equation}
\medskip

Let $E(t)$ denote the energy of $w$,
\begin{equation}\label{3.6}
E(t) = \frac{1}{2} \int |\nabla w|^{2} + \frac{1}{4} \int |w|^{4}.
\end{equation}
First observe that Lemma $\ref{l2.4}$ and $\| w \|_{\dot{H}^{1/2}} \lesssim 1$ implies that $E(1) < \infty$. To prove Theorem $\ref{t1.1}$, it suffices to prove that for any $T_{0} > 1$ such that $(\ref{1.1})$ has a solution on $[0, T_{0})$
\begin{equation}\label{3.6.1}
\sup_{t \in [1, T_{0})} E(t) < \infty.
\end{equation}
Indeed, by interpolation and the Sobolev embedding theorem, $E(t) < \infty$ implies that $\| w(t) \|_{L^{5}} < \infty$. Meanwhile, by $(\ref{2.15})$--$(\ref{2.17})$, $(\ref{2.11.1})$, and $(\ref{2.25})$, $\| \tilde{v}(t) \|_{L^{5}}$ is uniformly bounded on $\mathbb{R}$. Therefore, $(\ref{3.6.1})$ implies
\begin{equation}\label{3.6.2}
\| u \|_{L_{t,x}^{5}([0, T_{0}) \times \mathbb{R}^{3})} < \infty.
\end{equation}

To estimate the growth of $E(t)$, compute the derivative in time of the energy. By $(\ref{3.5})$,
\begin{equation}\label{3.7}
\frac{d}{dt} E(t) = -\langle \Delta w, w_{t} \rangle + \langle |w|^{2} w, w_{t} \rangle = \langle |w|^{2} w - |u|^{2} u, w_{t} \rangle,
\end{equation}
where $\langle \cdot, \cdot, \rangle$ is the inner product
\begin{equation}\label{3.8}
\langle f, g \rangle = Re \int f(x) \bar{g}(x) dx.
\end{equation}
By the product rule,
\begin{equation}\label{3.9}
\aligned
\langle w_{t}, |u|^{2} u - |w|^{2} w \rangle = \frac{d}{dt} \langle |w|^{2} w, \tilde{v} \rangle + \frac{d}{dt} \langle |\tilde{v}|^{2}, |w|^{2} \rangle + \frac{1}{2} \frac{d}{dt} Re \int \bar{w}^{2} \tilde{v}^{2} + \frac{d}{dt} \langle w, |\tilde{v}|^{2} \tilde{v} \rangle \\ - 2 \langle \tilde{v}_{t} \bar{\tilde{v}}, |w|^{2} \rangle - \langle |w|^{2} w, \tilde{v}_{t} \rangle - Re \int w^{2} \bar{\tilde{v}} \bar{\tilde{v}}_{t} - 2 \langle w, |\tilde{v}|^{2} \tilde{v}_{t} \rangle - \langle w, \tilde{v}^{2} \bar{\tilde{v}}_{t} \rangle.
\endaligned
\end{equation}
Then define the modified energy,
\begin{equation}\label{3.10}
\mathcal E(t) = E(t) - \langle |w|^{2} w, \tilde{v} \rangle - \langle |\tilde{v}|^{2}, |w|^{2} \rangle - \frac{1}{2} Re \int w^{2} \bar{\tilde{v}}^{2} - \langle w, |\tilde{v}|^{2} \tilde{v} \rangle.
\end{equation}
By Holder's inequality, and the fact that $\| \tilde{v} \|_{L^{4}} \lesssim_{\delta} 1$ for all $t \in [1, \infty)$ (again using $(\ref{2.15})$--$(\ref{2.17})$, $(\ref{2.11.1})$, and $(\ref{2.25})$),
\begin{equation}\label{3.11}
\langle |w|^{2} w, \tilde{v} \rangle + \langle |\tilde{v}|^{2}, |w|^{2} \rangle + \frac{1}{2} Re \int w^{2} \bar{\tilde{v}}^{2} + \langle w, |\tilde{v}|^{2} \tilde{v} \rangle \lesssim E(t)^{3/4} + E(t)^{1/4}.
\end{equation}
Therefore, when $E(t)$ is large, $E(t) \sim \mathcal E(t)$. Since we are attempting to prove a uniform bound for $E(t)$, it is enough to uniformly bound $\mathcal E(t)$.\medskip

Also, by $(\ref{3.9})$,
\begin{equation}\label{3.12}
\frac{d}{dt} \mathcal E(t) = -\langle |w|^{2} w, \tilde{v}_{t} \rangle - 2 \langle \tilde{v}_{t} \bar{\tilde{v}}, |w|^{2} \rangle - Re \int w^{2} \bar{\tilde{v}} \bar{\tilde{v}}_{t} - 2 \langle w, |\tilde{v}|^{2} \tilde{v}_{t} \rangle - \langle w, \tilde{v}^{2} \bar{\tilde{v}}_{t} \rangle.
\end{equation}
Since $\tilde{v}$ solves $(\ref{3.4})$, $\tilde{v}_{t} = i \Delta \tilde{v} = i \Delta u_{l} + i \Delta v$.\medskip

Lemma $\ref{l2.2}$ implies that for any $t > 1$,
\begin{equation}\label{3.12.2}
\| v(t) \|_{L^{\infty}} + \| \nabla v(t) \|_{L^{\infty}} = \| \nabla \int_{0}^{(1 - \delta)} e^{i(t - \tau) \Delta} \langle \nabla \rangle F(u) d\tau \|_{L^{\infty}} \lesssim \frac{1}{\delta^{3/2} t^{3/2}}.
\end{equation}
Therefore,
\begin{equation}\label{3.13}
\langle |w|^{2} w, v_{t} \rangle = \langle |w|^{2} w, i \Delta v \rangle = -\langle \nabla (|w|^{2} w), i \nabla v \rangle \lesssim \| \nabla v \|_{L^{\infty}} \| \nabla w \|_{L^{2}} \| w \|_{L^{4}}^{2} \lesssim_{\delta} \frac{1}{t^{3/2}} E(t).
\end{equation}
\noindent \textbf{Remark:} Since $\delta > 0$ is fixed, we will ignore it from now on.\medskip

\noindent Also, by Holder's inequality,
\begin{equation}\label{3.14}
\langle i \Delta (e^{it \Delta} u_{0}), |w|^{2} w \rangle \lesssim \| |\nabla|^{11/7} u_{l} \|_{L^{7}} \| \nabla w \|_{L^{2}}^{3/7} \| w \|_{L^{4}}^{18/7} \lesssim \frac{1}{t^{15/14}} E(t)^{6/7}.
\end{equation}
This takes care of the contribution of $\langle \tilde{v}_{t}, |w|^{2} w \rangle$.\medskip

Next, integrating by parts,
\begin{equation}\label{3.15}
2 \langle i (\Delta \tilde{v}) \bar{\tilde{v}}, |w|^{2} \rangle = -2 \langle i |\nabla \tilde{v}|^{2}, |w|^{2} \rangle - 2 \langle i (\nabla \tilde{v}) \bar{\tilde{v}}, \nabla |w|^{2} \rangle = -2 \langle i (\nabla \tilde{v}) \bar{\tilde{v}}, \nabla |w|^{2} \rangle.
\end{equation}
Then by Holder's inequality, since $\| \tilde{v} \|_{L^{4}} \lesssim 1$,
\begin{equation}\label{3.16}
\langle i (\nabla v) \bar{\tilde{v}}, \nabla |w|^{2} \rangle \lesssim \| \nabla v \|_{L^{\infty}} \| \tilde{v} \|_{L^{4}} \| w \|_{L^{4}} \| \nabla w \|_{L^{2}} \lesssim \frac{1}{t^{3/2}} E(t)^{3/4}.
\end{equation}
Also, by Holder's inequality and interpolation,
\begin{equation}\label{3.17}
\langle i (\nabla u_{l}) (u_{l}), \nabla |w|^{2} \rangle \lesssim \| \nabla u_{l} \|_{L_{x}^{\infty}} \| u_{l} \|_{L^{4}} \| \nabla w \|_{L^{2}} \| w \|_{L^{4}} \lesssim \frac{1}{t} \frac{1}{t^{1/8}} E(t)^{3/4}.
\end{equation}
Finally,
\begin{equation}\label{3.17.1}
\langle i (\nabla u_{l}) v, \nabla |w|^{2} \rangle \lesssim \| \nabla u_{l} \|_{L_{x}^{\infty}} \| v \|_{L^{3}}^{3/4} \| v \|_{L^{\infty}}^{1/4} \| \nabla w \|_{L^{2}} \| w \|_{L^{4}} \lesssim \frac{1}{t} \frac{1}{t^{3/8}} E(t)^{3/4}.
\end{equation}
In $(\ref{3.17})$ and $(\ref{3.17.1})$ we used
\begin{lemma}\label{l3.1}
\begin{equation}\label{3.21}
\| u_{l} \|_{L^{4}} \lesssim \frac{1}{t^{1/8}},
\end{equation}
and
\begin{equation}\label{3.22}
\| \nabla u_{l} \|_{L^{\infty}} \lesssim \frac{1}{t}.
\end{equation}
\end{lemma}
\noindent \emph{Proof:} This is proved by interpolating $(\ref{2.15})$--$(\ref{2.17})$. By Bernstein's inequality, $(\ref{2.16})$, $(\ref{2.17})$, and the Sobolev embedding theorem,
\begin{equation}
\| \nabla P_{\leq t^{-1/2}} u_{l} \|_{L^{\infty}} + \| \nabla P_{\geq t^{-1/2}} u_{l} \|_{L^{\infty}} \lesssim \frac{1}{t}.
\end{equation}
Also by the Bernstein inequality and the Sobolev embedding theorem, along with $(\ref{2.17})$ and $u_{l} \in \dot{H}^{1/2}$,
\begin{equation}
\| P_{\geq t^{-1/2}} u_{l} \|_{L^{4}} + \| P_{\leq t^{-1/2}} u_{l} \|_{L^{4}} \lesssim \frac{1}{t^{1/8}}.
\end{equation}
This proves the Lemma. $\Box$\medskip

The contribution of $2 Re \int w^{2} \bar{v} \bar{v}_{t}$ may be estimated in a similar manner as the contribution of $(\ref{3.15})$, except that there is an additional term to consider,
\begin{equation}\label{3.18}
-2 Re \int i w^{2} (\nabla \bar{\tilde{v}})^{2}.
\end{equation}
Interpolating $(\ref{3.22})$ with $(\ref{2.17})$,
\begin{equation}\label{3.20}
-2 Re \int i w^{2} (\nabla \bar{u}_{l})^{2} \lesssim \| \nabla u_{l} \|_{L^{4}}^{2} \| w \|_{L^{4}}^{2} \lesssim \frac{1}{t^{5/4}} E(t)^{1/2}.
\end{equation}
Meanwhile, replacing $\| |\nabla|^{1/2} P_{j} u \|_{L_{t}^{\infty} L_{x}^{2}}$ with $\| |\nabla|^{1/2} P_{j} u \|_{L_{t}^{2} L_{x}^{6}}$ in $(\ref{2.20})$--$(\ref{2.22})$ implies that for $t > 1$,
\begin{equation}
\| \nabla \int_{0}^{(1 - \delta)} e^{i(t - \tau) \Delta} F(u) \|_{L_{x}^{3}} \lesssim \frac{1}{ t^{1/2}}.
\end{equation}
Interpolating $(\ref{3.12.2})$ with $(\ref{3.21})$,
\begin{equation}\label{3.19}
-2 Re \int i w^{2} (\nabla \bar{\tilde{v}})^{2} \lesssim \| \nabla \tilde{v} \|_{L^{4}}^{2} \| w \|_{L^{4}}^{2} \lesssim \frac{1}{t^{3/2}} E(t)^{1/2}.
\end{equation}

Now treat
\begin{equation}\label{3.23}
2 \langle w, |\tilde{v}|^{2} \tilde{v}_{t} \rangle + \langle w, \tilde{v}^{2} \bar{\tilde{v}}_{t} \rangle = 2 \langle w, |\tilde{v}|^{2} (i \Delta \tilde{v}) \rangle + \langle w, \tilde{v}^{2} \overline{(i \Delta \tilde{v})} \rangle.
 \end{equation}
 After integrating by parts,
 \begin{equation}\label{3.24}
 \aligned
 (\ref{3.23}) \lesssim \langle |\nabla \tilde{v}|^{2}, |v| |w| \rangle + \langle |\nabla \tilde{v}| |\nabla w|, |v|^{2} \rangle \\ \lesssim \| \nabla \tilde{v} \|_{L^{4}}^{2} \| \tilde{v} \|_{L^{4}} \| w \|_{L^{4}} + \| \nabla w \|_{L^{2}} \| \nabla \tilde{v} \|_{L^{\infty}} \| \tilde{v} \|_{L^{4}}^{2} \lesssim \frac{1}{t^{5/4}} E(t)^{1/4} + \frac{1}{t} E(t)^{1/2} \| \tilde{v}(t) \|_{L^{4}}.
 \endaligned
 \end{equation}
 Interpolating $(\ref{3.12})$ with $\| v \|_{L^{3}} \lesssim 1$ implies $\| v \|_{L^{4}} \lesssim t^{-3/8}$, which combined with $(\ref{3.21})$ implies $\| \tilde{v} \|_{L^{4}} \lesssim \frac{1}{t^{1/8}}$. Therefore, we have proved
 \begin{equation}\label{3.25}
 \frac{d}{dt} \mathcal E(t) \lesssim \frac{1}{t^{15/14}} (1 + \mathcal E(t)).
 \end{equation}
By Gronwall's inequality, $(\ref{3.25})$ implies a uniform bound on $\mathcal E(t)$. This implies a uniform bound on $E(t)$, since $E(t) \sim \mathcal E(t)$ when $E(t)$ is large, which proves Theorem $\ref{t1.1}$.

 \section*{Acknowledgments}
The author was supported by NSF grants DMS-1500424 and DMS-1764358 during the time of writing this paper. The author is also grateful to Frank Merle and Kenji Nakanishi for helpful conversations regarding the defocusing nonlinear Schr{\"o}dinger equation.

\bibliography{biblio}
\bibliographystyle{alpha}
\medskip

\end{document}